\documentclass[12pt]{amsart}
\usepackage{amssymb}

{\theoremstyle{definition}
    \newtheorem{theorem}{Theorem}[section]
    \newtheorem{proposition}[theorem]{Proposition}
    
    \newtheorem{corollary}[theorem]{Corollary}
    \newtheorem{definition}[theorem]{Definition}
    \newtheorem{remark}[theorem]{Remark}
    
    \newtheorem{parrafo}[theorem]{{\!}}
}

\numberwithin{equation}{theorem}

\newcommand{\calo}{{\mathcal {O}}}

\newcommand{\weak}{{\curlyvee}}

\DeclareMathOperator{\codim}{codim}
\DeclareMathOperator{\Diff}{Diff}
\DeclareMathOperator{\Sing}{Sing}
\DeclareMathOperator{\MaxB}{\underline{\mathbf{Max}}}
\DeclareMathOperator{\ord}{ord}
\DeclareMathOperator{\Tay}{Tay}
\DeclareMathOperator{\word}{w-ord}

\title[]{Rees algebras and resolution of singularities}

\author{Santiago Encinas}
\author{Orlando Villamayor}

\address{Dpto. Matem\'atica Aplicada. Universidad de Valladolid.
47014 Valladolid. Spain.}
\email{sencinas@maf.uva.es}
\address{Dpto. Matem\'aticas, Facultad de Ciencias, Universidad
Aut\'onoma de Madrid, Canto Blanco 28049 Madrid, Spain.}
\email{villamayor@uam.es}
\thanks{2000 {\em Mathematics subject classification. 14E15.}}
\thanks{.}

\subjclass{}

\keywords{Resolution of singularities. Desingularization}
\date{August 2006} \dedicatory{ } \commby{}

\begin{document}
\begin{abstract}

Embedded principalization of ideals in smooth schemes, also known as
Log-resolutions of ideals, play a central role in algebraic geometry.
If two sheaves of ideals, say $I_1$ and $I_2$, over a smooth scheme
$V$ have the same integral closure, it is well known that
Log-resolution of one of them induces a Log-resolution of the other.
On the other hand, in case $V$ is smooth over a field of
characteristic zero, an algorithm of desingularization provides, for
each sheaf of ideals, a unique Log-resolution.

In this paper we show that algorithms of desingularization define the
same Log-resolution for two ideals having the same integral closure.
We prove this result here by using the form of induction introduced by
W{\l}odarczyk.

We extend the notion of Log-resolution of ideals over a smooth scheme
$V$, to that of Rees algebras over $V$; and then we show that two Rees
algebras with the same integral closure undergo the same constructive
resolution.  The key point is the interplay of integral closure with
differential operators.

\end{abstract}
\maketitle


\tableofcontents


\section*{Introduction.}
\label{introduction}  Differential operators on smooth schemes have
played a central role in the study of embedded desingularization.

J. Giraud provides an alternative approach to the form of
induction used by Hironaka in his Desingularization Theorem (over
fields of characteristic zero).  In doing so, Giraud introduced
technics based on differential operators (\cite{Gir},
\cite{Giraud1975}). This result was important for the development
of algorithms of desingularization in the late 80's (i.e. for
constructive proofs of Hironaka's theorem).

Differential operators appeared in the work of J. W{\l}odarczyk
(\cite{Wlodarczyk2005}), and also on the notes of J. Koll\'{a}r
(\cite{Kollar2005prep}); where algorithms of resolution are
developed.

The notions of Rees algebras over smooth schemes, and that of Rees
algebras closed by higher order differentials, already appear in Hironakas
study on infinitely near points (\cite{Hironaka03};
\cite{Hironaka05} ), and more recently in Kawanoue's work in \cite{Kawanoue}.

A Log-resolution, or embedded principalization, of an ideal $I$ on
a smooth scheme $V$, is a proper birational morphism of smooth
schemes, say $V'\to V$, so that the total transform of $I$ is an
invertible ideals in $V'$ supported on smooth hypersurfaces having
only normal crossings. When $V$ is smooth over a field of
characteristic zero there are algorithms that provide a
Log-resolution of an ideal $I$. We shall make use of Rees algebras in
proving that two ideals with the same integral closure undergo the
same algorithmic Log-resolution (\ref{teofin}).

The paper is organized so as to motivate the extension of
Log-resolution theorems of ideals over fields of characteristic
zero, to the case Rees algebras, this is done in Sections 1 and 2.

In Sections 3 and 4, the reader is introduced to the fascinating
relation of differential operators acting on Rees algebras, with
the notion of integral closure of these algebras. These first 4 sections
are included for self-containment. We refer to
\cite{Villamayor2006prep1}, or \cite{Villamayor2006prep3}, for details.

 In Section 5 we discuss some natural equivalence relation on Rees
algebras when it comes to desingularization. Finally, in sections
6 and 7 we discuss the main results.

In this paper we always consider smooth schemes over fields of
characteristic zero, however the extension of resolution theorems to Rees algebras, treated in this work, is also motivated by recent
development of invariants over arbitrary fields.

We refer to
\cite{Villamayor2006prep2} where a link of differential operators with elimination
theory is presented. In that paper elimination of one variable is formulated in terms of Rees algebras.
Over fields of characteristic zero this elimination recovers  Hironaka«s form of induction in desingularization theorems. However new invariants arise from this form of elimination, defined entirely in terms of Rees algebras, over fields of positive characteristic.


\section{Monoidal transformations and Hironaka's topology.}
\label{sec1}

\begin{parrafo}\label{parA}
Fix a smooth scheme \(V\) over a field \(k\), an ideal \(J\subset
\calo_V\), and a positive integer \(b\).  Hironaka attaches to these
data, say \((J,b)\), a Zariski closed set in $V$, say
\begin{equation*}
 \Sing(J,b):=\{x \in V/ \nu_x(J_x)\geq b\}
\end{equation*}
where \(\nu_x(J_x)\) denotes the order of \(J\) at the local regular
ring \(\calo_{V,x}\).

Given \((J_{1},b_{1})\) and \((J_{2},b_{2})\), then
\begin{equation*}
\Sing(J_{1},b_{1})\cap\Sing(J_{2},b_{2})=\Sing(K,c)
\end{equation*}
where \(K=J_{1}^{b_{2}}+J_{2}^{b_{1}}\), and \(c=b_{1}\cdot b_{2}\).
Set formally \((J_{1},b_{1})\odot (J_{2},b_{2})=(K,c)\).

There is also a notion of {\em permissible transformation} on these
data \((J,b)\).  Let \(Y\) be a smooth subscheme in \(V\), included in
the closed \(\Sing(J,b)\), and let
\begin{equation}\label{eq1}
\begin{array}{ccl}
V & \stackrel{\pi }{\longleftarrow}& V'
\\
 \cup & & \cup \\
  Y & & H'=\pi^{-1}(Y) ,\\
\end{array}
\end{equation}
be the blow up of \(V\) at a smooth sub-scheme \(Y\).  Note that
\begin{equation*}
J\calo_{V'}=I(H')^b J',
\end{equation*}
where \(I(H')\) is the sheaf of functions vanishing along the
exceptional hypersurface \(H'\).

We call \((J',b)\) the {\em transform} of \((J,b)\) by the permissible
monoidal transformation.

If \(\pi\) is permissible for both \((J_{1},b_{1})\) and
\((J_{2},b_{2})\), then it is permissible for \((K,c)\).  Moreover, if
\((J'_1,b_{1})\), \((J'_2,b_{2})\), and \((K',c)\) denote the
transforms, then \((J'_1,b_{1})\odot (J'_2,b_{2})=(K',c)\).
\end{parrafo}

\begin{parrafo}\label{parB}
We will consider $\mathbb{N}$-graded algebras. Fix a variable $W$ and define a {\em Rees algebra} over \(V\) to be a graded
noetherian subring of \(\calo_V[W]\), say:
\begin{equation*}
\mathcal{G}=\bigoplus_{n\geq 0}I_{n}W^n,
\end{equation*}
where \(I_0=\calo_V\) and each \(I_n\) is a sheaf of ideals.  By assumption, at every affine open set \(U\subset V\) there is a finite
set
\begin{equation*}
\mathcal{F}=\{f_1W^{n_1},\dots ,f_sW^{n_s}\}, \end{equation*}
\(n_i\geq 1\) and \(f_i\in \calo_V(U)\), so that the restriction of
\(\mathcal{G}\) to \(U\) is
\begin{equation*}
\calo_V(U)[f_1W^{n_1},\dots ,f_sW^{n_s}] (\subset \calo_V(U)[W]).
\end{equation*}
To a Rees algebra \(\mathcal{G}\) we attach a closed set:
\begin{equation*}
    \Sing(\mathcal{G}):=\{ x\in V/ \nu_x(I_n)\geq n, \text{ for every }
    n\geq 1\},
\end{equation*}
where \(\nu_x(I_n)\) denotes the order of the ideal \(I_n\) at the
local regular ring \(\calo_{V,x}\).
\end{parrafo}

\begin{remark}\label{rkK1}
Rees algebras are related to Rees rings of ideals.  A Rees algebra
is a Rees ring if, given any affine open set \(U\subset V\), and
\(\mathcal{F}=\{f_1W^{n_1},\dots ,f_sW^{n_s}\}\) as above, all
degrees \(n_i\) are one.  In such case it is the Rees ring of the
ideal \(I=\langle f_1,\dots ,f_s \rangle\).  The integral closure
of a Rees ring of an ideal is no longer a Rees ring of another
ideal, however it is within the class of rings we consider here:
the integral closure of a Rees ring is a Rees algebra.

In general Rees algebras are, in some sense, integral over Rees rings.
In fact, if \(N\) is a positive integer divisible by all \(n_i\), it
is easy to check that
\begin{equation*}
\calo_V(U)[f_1W^{n_1},\dots,f_sW^{n_s}]= \bigoplus_{n\geq 0}I_{n}W^n
(\subset \calo_V(U)[W]), \end{equation*} is integral over the Rees
sub-ring \(\calo_V(U)[I_NW^N](\subset \calo_V(U)[W^N])\).  In fact,
\(\calo_V(U)[I_NW^N]\subset \calo_V(U)[f_1W^{n_1},\dots
,f_sW^{n_s}]\), and \((f_iW^{n_i})^{\frac{N}{n_i}}\in
\calo_V(U)[I_NW^N]\).
\end{remark}

\begin{proposition}\label{prop1}
Given an affine open \(U\subset V\), and
\(\mathcal{F}=\{f_1W^{n_1},\dots ,f_sW^{n_s}\}\) as above,
\begin{equation*}
\Sing(\mathcal{G})\cap U= \bigcap_{1\leq i \leq s}\{x\in U\mid
\nu_{x}(f_i)\geq n_i\}.
\end{equation*}

\end{proposition}
\begin{proof}
It is clear that \(\nu_x(f_i)\geq n_i\) for \(x\in
\Sing(\mathcal{G})\), \(0\leq i \leq s\).  So
\begin{equation*}
\Sing(\mathcal{G})\cap U\subset \bigcap_{1\leq i \leq s}\{x\in U\mid
\nu_{x}(f_i) \geq n_i\}.
\end{equation*}
On the other hand, for every index \(N\geq 1\), \(I_N(U)W^N\) is
generated by elements of the form \(G_N(f_1W^{n_1},\dots
,f_sW^{n_s})\), where \(G_N(Y_1,\dots ,Y_s)\in \calo_U[Y_1,\dots
,Y_s]\) is weighted homogeneous of degree \(N\), provided each \(Y_j\)
has weight \(n_j\).  The reverse inclusion is now clear.
\end{proof}

\begin{parrafo}\label{parwt}
A monoidal transformation (\ref{eq1}) is said to be {\em permissible}
for \(\mathcal{G}\) if \(Y\subset\Sing(\mathcal{G})\).  In such case,
for each index \(n\geq 1\), there is a sheaf of ideals, say
\(I'_{n}\subset \calo_{V'}\), so that
\begin{equation*}
I_n\calo_{V'}= I(H')^n I'_n.
\end{equation*}

We define the {\em total transform} of \(\mathcal{G}\) to be
\(\bigoplus_{n\geq 0}I_n\calo_{V'}W^n\).  On the other hand we define
{\em weighted transform} of \(\mathcal{G}\) as:
\begin{equation*}
\mathcal{G}'=\bigoplus_{n\geq 0}I'_{n}W^n;
\end{equation*}
which is a Rees algebra over \(V'\) (see \ref{prop2}).

Let \(\mathcal{G}=\bigoplus_{n\geq 0}I_{n}W^n\) be a Rees algebra on
\(V\), and set \(V\stackrel{\pi }{\longleftarrow} V'\) a permissible
transformation of \(\mathcal{G}\).  Let \(U\subset V\) be affine open
set, and \(\mathcal{F}=\{f_1W^{n_1},\dots ,f_sW^{n_s}\}\) be such that
the restriction of \(\mathcal{G}\) to \(U\) is \(
\calo_V(U)[f_1W^{n_1},\dots ,f_sW^{n_s}] (\subset \calo_V(U)[W]).  \)
Note that the total transform of \(\mathcal{G}\), restricted to the
open set \(\pi^{-1}(U)(\subset V')\), is also generated by
\(\{f_1W^{n_1},\dots ,f_sW^{n_s}\}(\subset
\calo_{V'}({\pi}^{-1}(U))[W])\).
\end{parrafo}

\begin{proposition}\label{prop2}
With the setting as above, there is an open covering of
\(\pi^{-1}(U)\) by affine sets \(U^{(\ell)}\), so that:
\begin{enumerate}
    \item \(\langle f_i \rangle =I(H'\cap U^{(l)})^{n_i}\langle
    f'_i\rangle \) for suitable \(f'_i\in \calo_{V'}(U^{(\ell)})\).

    \item The restriction of the weighted transform, say
    \(\mathcal{G}'\), to each open set \(U^{(\ell)}\) is
    \begin{equation*}
    \calo_{V'}(U^{(\ell)})[f'_1W^{n_1},\dots ,f'_sW^{n_s}] (\subset
    \calo_{V'}(U^{(\ell)})[W]).
    \end{equation*}
\end{enumerate}
\end{proposition}

\begin{proof}
(1) Follows from proposition~\ref{prop1}, since every \(f_{i}\) has
order at least \(n_i\) along the center \(Y\).  For (2) argue as
in the proof of Proposition~\ref{prop1}, by using the fact that
each ideal \(I_N\) is generated by weighted homogeneous
polynomials on the element of \(\mathcal{F}\).
\end{proof}

Given two Rees algebras over \(V\), say
\(\mathcal{G}_{1}=\bigoplus_{n\geq 0}I_{n}W^n\) and
\(\mathcal{G}_{2}=\bigoplus_{n\geq 0}J_{n}W^n\), set \(K_n=I_n+J_n\)
in \(\calo_V\), and define:
\begin{equation*}
\mathcal{G}_{1}\odot\mathcal{G}_{2}=\bigoplus_{n\geq 0}K_{n}W^n,
\end{equation*}
as the subalgebra of \(\calo_V[W]\) generated by \(\{K_{n}W^n, n\geq
0\}\).

One can check that:
\begin{enumerate}
    \item \(\Sing(\mathcal{G}_{1}\odot\mathcal{G}_{2})=
    \Sing(\mathcal{G}_{1})\cap\Sing(\mathcal{G}_{2})\).  In
    particular, if \(\pi\) in (\ref{eq1}) is permissible for
    \(\mathcal{G}_{1}\odot\mathcal{G}_{2}\), it is also permissible
    for \(\mathcal{G}_{1}\) and for \( \mathcal{G}_{2}\).

    \item Set \(\pi\) as in (1), and let
    \((\mathcal{G}_{1}\odot\mathcal{G}_{2})'\), \(\mathcal{G}'_1\),
    and \(\mathcal{G}'_2\) denote the transforms at \(V'\).  Then:
    \begin{equation*}
    (\mathcal{G}_{1}\odot\mathcal{G}_{2})'=
    \mathcal{G}'_1\odot\mathcal{G}'_2.
    \end{equation*}
\end{enumerate}

\section{On Hironaka pairs and Rees algebras.}
\label{sec2}

Recall that two ideals, say \(I\) and \(J\), in a normal domain \(R\)
have the same integral closure if they are equal for every extension to
a valuation ring (i.e. if \(IS=JS\) for every ring homomorphism \(R\to
S\) on a valuation ring \(S\)).  \medskip

Hironaka considers the following equivalence on pairs \((J,b)\) over a
smooth scheme \(V\).
\begin{definition}\label{defe}
The pairs \((J_{1},b_{1})\) and \((J_{2},b_{2})\) are {\em idealistic}
equivalent on \(V\) if \(J_{1}^{b_{2}}\) and \(J_{2}^{b_{1}}\) have
the same integral closure.
\end{definition}

Among Rees algebras the equivalence relation, also defined in terms of
integral closure, is:
\begin{definition}\label{defei}
We say that two Rees algebras over \(V\), say
\(\mathcal{G}_{1}=\bigoplus_{n\geq 0}I_{n}W^n\) and
\(\mathcal{G}_{2}=\bigoplus_{n\geq 0}J_{n}W^n\), are {\em integrally
equivalent}, if both have the same integral closure in \(\calo_V[W]\).
\end{definition}
In general we want to identify two Rees algebras if they have the same
integral closure.  This notion will be revisited in
section~\ref{sec5}, where it will be linked with a weaker equivalence
relation.

\begin{parrafo}\label{parC}
We assign to a pair \((J,b)\) over a smooth scheme \(V\) the Rees
algebra, say:
\begin{equation*}
\mathcal{G}_{(J,b)}=\calo_V[JW^b].
\end{equation*}
Note that \(\mathcal{G}_{(J,b)}\) is a Rees ring of an ideal in
\(\calo_V[W^b]\), but we can consider it as a graded subalgebra in
\(\calo_V[W]\).  Remark \ref{rkK1} shows that every Rees algebra is, in
this sense, integrally equivalent to the Rees ring attached to a pair.
In fact if \(\mathcal{G}=\bigoplus_{n\geq 0}I_{n}W^n\), then it has
the same integral closure as \(\mathcal{G}_{(I_N,N)}\) for a suitable
\(N\).
\end{parrafo}

\begin{parrafo} \label{Require}

A key point in our development is to attach invariants or
geometric objects to Rees algebras. This will be always be done
subject to the following two requirements:

\begin{enumerate}
    \item Every construction or invariant attached to a Rees algebra will be the same for
    two integrally equivalent Rees algebras.

    \item To all construction and invariants we present for Rees algebras, there will be a similar one on the
    class of idealistic pairs.
\end{enumerate}
\end{parrafo}

For example the operator \( \odot \) fulfills our requirement:
\begin{proposition}\label{prop3}
Set \(\mathcal{G}_{1}=\mathcal{G}_{(J_{1},b_{1})}\) and
\(\mathcal{G}_{2}=\mathcal{G}_{(J_{2},b_{2})}\) (i.e. the Rees
algebras corresponding to Hironaka's pairs \((J_{1},b_{1})\) and
\((J_{2},b_{2})\)), then \(\mathcal{G}_{1}\odot\mathcal{G}_{2}\) is
equivalent to the Rees algebra assigned to
\((J_{1},b_{1})\odot(J_{2},b_{2})\).  Furthermore, this relation is
preserved by transformations.
\end{proposition}

\begin{proof}
Fix an affine open set \(U\) in \(V\), \(\{f_1,\dots , f_s\}\in
\calo_V(U)\) generators of \(J_{1}(U)\), and \(\{g_1,\dots , g_r\}\in
\calo_V(U)\) generators of \(J_{2}(U)\).  Then:
\begin{description}
    \item[(i)] The restriction of \(\mathcal{G}_{1}\) to \(U\) is
    \begin{equation*}
    \calo_V(U)[f_1W^{b_{1}},\ldots,f_sW^{b_{1}}] (\subset
    \calo_V(U)[W]).
    \end{equation*}

    \item[(ii)] The restriction of \(\mathcal{G}'\) is
    \begin{equation*}
    \calo_V(U)[g_1W^{b_{2}},\ldots,g_rW^{b_{2}}] (\subset
    \calo_V(U)[W]).
    \end{equation*}

    \item[(iii)] The restriction of
    \(\mathcal{G}_{1}\odot\mathcal{G}_{2}\) to \(U\) is
    \begin{equation*}
    \calo_V(U)[ f_1W^{b_{1}},\ldots,f_sW^{b_{1}},
    g_1W^{b_{2}},\ldots,g_rW^{b_{2}}] (\subset \calo_V(U)[W]).
    \end{equation*}

    \item[(iv)] The restriction of the Rees algebra assigned to
    \((J_{1},b_{1})\odot(J_{2},b_{2})\) is generated by
    \begin{equation*}
    \{ (f_1^{\alpha_1}\cdots f_s^{\alpha_s})\cdot W^{b_{1}b_{2}};
    (g_1^{\beta_1}\cdots g_s^{\beta_s})\cdot W^{b_{1}b_{2}} \mid
    \alpha_1+\cdots+\alpha_s=b_{2};\ \beta_1+\cdots +\beta_r=b_{1}\}.
    \end{equation*}
\end{description}

One can finally check that both algebras in (iii) and (iv) have the
same integral closure in \(\calo_V(U)[W]\).  The last assertion,
namely that the relation is preserved by transformations, is now
straight forwards.
\end{proof}

\begin{proposition}
Let \( \mathcal{G}_{ij} \), \( i,j=1,2 \), be Rees algebras such that
\( \mathcal{G}_{1j} \) and \( \mathcal{G}_{2j} \) are integrally
equivalent, \( j=1,2 \). Then the Rees algebras \(
\mathcal{G}_{11}\odot\mathcal{G}_{12} \) and \(
\mathcal{G}_{21}\odot\mathcal{G}_{22} \) are integrally equivalent.
\end{proposition}

\begin{definition}\label{dver}
Let \( \mathcal{G}=\bigoplus_{n}I_{n}W^{n} \) be a Rees algebra.  For
every \( m\in\mathbb{Z} \), \( m>0 \) we set
\begin{equation}\label{4eq1}
    V^{(m)}(\bigoplus_{n}I_{n}W^{n})=\bigoplus_{n\geq0} I_{mn}W^{mn}
    (\subset\calo_{V}[W^m])
\end{equation}
\end{definition}

\begin{parrafo}\label{paramigo}
For every \(f_nW^n \in \mathcal{G}\), \((f_nW^n)^m\in
V^{(m)}(\mathcal{G})\); which shows that \(
V^{(m)}(\mathcal{G})\subset\mathcal{G} \) is an integral extension for
all \( m \).  So that \( \mathcal{G} \) and \( V^{(m)}(\mathcal{G}) \)
are integrally equivalent (\ref{defei}).

The following properties hold for \(V^{(m)}\):
\begin{enumerate}
    \item If \( \mathcal{G}_{1} \) and \( \mathcal{G}_{2} \) are
    integrally equivalent Rees algebras then \(
    V^{(m)}(\mathcal{G}_{1}) \) and \( V^{(m)}(\mathcal{G}_{2}) \) are
    integrally equivalent (i.e. \( V^{(m)}\) is compatible with
    integral closure (\ref{Require})).

    \item If \( \mathcal{G}=\mathcal{G}_{(J,b)} \) then \(
    V^{(m)}(\mathcal{G})$ and $\mathcal{G}_{(J^{m},bm)} \) are integrally equivalent.
\end{enumerate}
In fact, if \(\bigoplus_{n}I_{n}W^{n}\subset
\bigoplus_{n}J_{n}W^{n}\) is the extension defined by taking
integral closures, then \(V^{(m)}(\bigoplus_{n}J_{n}W^{n})\) is
the integral closure of \(V^{(m)}(\bigoplus_{n}I_{n}W^{n})\) in
\(\calo_{V}[W^m]\).

In the particular case in which the Rees algebra is the Rees ring of
an ideal, it turns out that \(J_n\) is the integral closure of
\(I_n\).  In general, each \(J_n\) contains the integral closure of
\(I_n\).

In \ref{rkK1} we have shown that for infinitely many suitable choices
of \( N \): \( V^{(N)}(\mathcal{G})\) is integral over the Rees ring
\(\calo_{V}[I_{N}W^{N}]\); namely, that \(\mathcal{G}_{(I_N,N)}
\subset \mathcal{G}\) is a (finite) integral extension.

Given a pair \( (J,b) \) and a positive integer \( m \), the pairs
\( (J,b) \) and \( (J^{m},bm) \) are idealistic equivalent
(\ref{defe}).  (1) and (2) show that \(\mathcal{G}_{(J,b)} \) and
\(\mathcal{G}_{(J^m,mb)}\) are integrally equivalent (\ref{defei}).

As all algebras considered here are finitely generated over \(\calo_V\), for \(N\)
suitably big: \( V^{(N)}(\mathcal{G})=\calo_{V}[I_{N}W^{N}]\).

\end{parrafo}
The next result follows from the previous discussion.
\begin{proposition} \label{PropEquivGJ}
Two pairs \( (J_{1},b_{1}) \) and \( (J_{2},b_{2}) \) are idealistic
equivalent if and only if the associated Rees algebras \(
\mathcal{G}_{1}=\mathcal{G}_{(J_{1},b_{1})} \) and \(
\mathcal{G}_{2}=\mathcal{G}_{(J_{2},b_{2})} \) are integrally
equivalent.
\end{proposition}

\begin{proof}
Note that the following are equivalent:
\begin{itemize}
    \item The pairs \( (J_{1},b_{1}) \) and \( (J_{2},b_{2}) \) are
    idealistic equivalent.

    \item The ideals \( J_{1}^{b_{1}} \) and \( J_{2}^{b_{1}} \) have
    the same integral closure.

    \item The Rees algebras \( V^{(b_{2})}(\mathcal{G}_{1}) \) and \(
    V^{(b_{1})}(\mathcal{G}_{2}) \) have the same integral closure.

    \item The Rees algebras \( \mathcal{G}_{1} \) and \(
    \mathcal{G}_{2} \) have the same integral closure.
\end{itemize}
\end{proof}

\section{On differential Rees algebras and Koll\'{a}r's tuned ideals.}
\label{sec3}

Here \(V\) is smooth over a field \(k\), so for each non-negative
integer \(r\) there is a locally free sheaf of differential operators
of order \(r\), say \(\Diff^{(r)}_k\).

\begin{definition} \label{3def1}
We say that a Rees algebra \(\mathcal{G}=\bigoplus I_nW^n\) is a
differential Rees algebra,  or simply a Diff-algebra, relative to
the field \(k\), if:
\begin{description}
    \item[(i)] \(I_n\supset I_{n+1}\).

    \item[(ii)] There is open covering of \(V\) by affine open sets
    \(\{U_i\}\), and for every \(D\in\Diff^{(r)}(U_i)\), and  \(h\in
    I_n(U_i)\), then \(D(h)\in I_{n-r}(U_i)\) provided \(n\geq r\).
\end{description}
\end{definition}

Given a sheaf of ideals \(I\subset\calo_V\) there is a natural
definition of an extension, say \(\Diff^{(r)}(I)\) (see Introduction).
Note that (ii) can be reformulated by
\begin{description}
    \item[(ii)'] \(\Diff^{(r)}(I_n)\subset I_{n-r}\) for each \(n\),
    and \(0\leq r \leq n\).
\end{description}

Diff-algebras are called differential structures in
\cite{Villamayor2006prep3}.

\begin{remark}
As we will view Rees algebras up to integral closure, it is not hard
to check that condition (i) can be imposed for Rees algebras.  In
fact, given \(\mathcal{G}=\bigoplus_{n\geq 0} I_n W^n \), we define
\(\mathcal{G}^{\natural}=\bigoplus I^{\natural}_n W^n \) by setting
\begin{equation*}
 I^{\natural}_n=\sum_{r\geq n} I_r.
\end{equation*}

If \(\mathcal{G}=\bigoplus I_n W^n \) is generated by \(\mathcal{F}=\{
g_{n_i}W^{n_i}, 1\leq i\leq m, n_i>0 \}\).  Namely, if:
\begin{equation*}
    \bigoplus_{n\geq 0} I_n\cdot W^n=
    \calo_{V}[\{g_{n_i}W^{n_i}\}_{i=1}^{m}],
\end{equation*}
then \(\mathcal{G}^{\natural}=\bigoplus I^{\natural}_n W^n \) is
generated by the finite set \(\{ g_{n_i}W^{n'_i}, 1\leq i\leq m, 1
\leq n'_i\leq n_i \}\),

Note that \(I^{\natural}_n \supset I^{\natural}_{n+1}\), and that
\(\bigoplus_{n} I_n W^n \subset \bigoplus I^{\natural}_n W^n \) is a
finite extension.  In fact, it suffices to check that given an element
\(g\in I_{n}\), then \(g W^{n-1}\) is integral over
\(\bigoplus_{n}I_n W^n \).  One can check that
\begin{equation*}
g\in I_{n}\Rightarrow g^{n-1} \in I_{n(n-1)}\Rightarrow g^{n} \in
I_{n(n-1)},
\end{equation*}
so \(g\cdot W^{n-1}\) fulfills the equation \(Z^{n}-(g^{n}\cdot
W^{n(n-1)})=0\).
\end{remark}

\begin{parrafo}\label{parod}
Fix a closed point \(x\in V\), and a regular system of parameters
\(\{x_1,\dots,x_d\}\) at \(\calo_{V,x}\).  The residue field, say
\(k'\) is a finite extension of \(k\), and the completion
\(\hat{\calo}_{V,x}=k'[[x_1,\dots,x_d]].\)

The Taylor development is the continuous \(k'\)-linear ring
homomorphism:
\begin{equation*}
\Tay: k'[[x_1,\dots,x_d]]\to k'[[x_1,\dots,x_d,T_1,\dots,T_d]]
\end{equation*}
that maps \(x_i\) to \(x_i+T_i\), \(1\leq i \leq d\).  So for
\(f\in k'[[x_1,\dots,x_d]]\),
\(\Tay(f(x))=\sum_{\alpha\in\mathbb{N}^d} g_{\alpha}T^{\alpha}\),
with \(g_{\alpha} \in k'[[x_1,\dots , x_d]]\). Define, for each
\(\alpha \in \mathbb{N}^d\), \(\Delta^{\alpha}(f)=g_{\alpha}\).
It turns out that
\begin{equation*}
\Delta^{\alpha}(\calo_{V,x})\subset \calo_{V,x},
\end{equation*}
and that \(\{ \Delta^{\alpha}, \alpha \in\mathbb{N}^d , 0 \leq
|\alpha| \leq r\}\) generate the \(\calo_{V,x}\)-module
\(\Diff^{(r)}_k(\calo_{V,x})\) (i.e. generate \(\Diff^{(r)}_k\)
locally at \(x\)).
\end{parrafo}

We finally introduce an operator on Rees algebras, which parallels
Giraud's extensions of ideals via differential operators.
\begin{theorem}\label{thopG}
For every Rees algebra \(\mathcal{G}\) over a smooth scheme \(V\),
there is a Diff-algebra, say \(G(\mathcal{G})\),  such that:
\begin{description}
    \item[(i)] \(\mathcal{G}\subset G(\mathcal{G})\).

    \item[(ii)] If \(\mathcal{G}\subset \mathcal{G}'\) and
    \(\mathcal{G}'\) is a Diff-algebra, then \(G(\mathcal{G})\subset
    \mathcal{G}'\).
\end{description}
Furthermore, if \(x\in V\) is a closed point, and
\(\{x_1,\dots,x_d\}\) is a regular system of parameters at
\(\calo_{V,x}\), and \(\mathcal{G}\) is locally generated by
\begin{equation*}
\mathcal{F}=\{ g_{i}W^{n_i}, n_i>0 , 1\leq i\leq m \},
\end{equation*}
then
\begin{equation}
\mathcal{F'}=\{\Delta^{\alpha}(g_{i})W^{n'_i-\alpha} \mid
g_{i}W^{n_i}\in \mathcal{F}, \alpha=(\alpha_1, \alpha_2,\dots ,
\alpha_d) \in\mathbb{N}^d, \quad 0\leq |\alpha| < n'_i \leq n_i\}
\end{equation}
generates \(G(\mathcal{G})\) locally at \(x\).
\end{theorem}

\begin{remark}
The local description in the Theorem shows that
\(\Sing(\mathcal{G})=\Sing(G(\mathcal{G}))\).

In fact, as \(\mathcal{G}\subset G(\mathcal{G})\), it is clear that
\(\Sing(\mathcal{G})\supset\Sing(G(\mathcal{G}))\).  For the converse
note that if \(\nu_x(g_{n_i})\geq n_i\), then
\(\Delta^{\alpha}(g_{n_i})\) has order at least \(n_i-|\alpha|\) at
the local ring \(\calo_{V,x}\).
\end{remark}

\begin{parrafo}\label{parKol}
In general \(\mathcal{G}\subset G(\mathcal{G})\), and equality
holds if \(\mathcal{G}\) is already a Diff-algebra.

Let \(\mathcal{G}=\bigoplus_{n\geq 0}I_nW^n\) be a Diff-algebra,
in particular it is integral over a Rees subring, say
\(\calo_V[I_NW^N]\) for suitable \(N\) (see \ref{rkK1}).  These
ideals \(I_N\) are called {\em tuned ideals} in
\cite{Kollar2005prep}, page 45.
\end{parrafo}

The previous Theorem defines an operator \(G\) that extends Rees
algebras into Diff-algebras.  Another natural operator we have
considered on Rees algebras it that defined by taking
normalization. The next Theorem relates both notions of
extensions.

\begin{theorem} \label{ThGEntera} \cite[Th 6.12]{Villamayor2006prep1}
Let \(\mathcal{G}_{1}\) and \(\mathcal{G}_{2}\) be integrally
equivalent Rees algebras on a smooth scheme \(V\), then
\(G(\mathcal{G}_{1})\) and \(G(\mathcal{G}_{2})\) are also integrally
equivalent (\ref{defei}).
\end{theorem}

\section{On differential Rees algebras and monoidal transformations.}
\label{sec4}

Let us briefly recall some previous results, where now \(J\subset
\calo_V\) be the sheaf of ideals defining a hypersurface \(X\) in the
smooth scheme \(V\).

So \(\Diff^0(J)=J\), and for each positive integer \(s\) there is an
inclusion \(\Diff^s(J)\subset\Diff^{s+1}(J)\) as sheaves of ideals in
\(\calo_V\), and hence \(V(\Diff^s(J))\supset V(\Diff^{s+1}(J))\).

Recall that \(b\) is the highest multiplicity at points of \(X\), if
and only if \(V(\Diff^{b}(J))=\emptyset\) and \(V(\Diff^{b-1}(J))\neq
\emptyset\) (i.e. if and only if \(\Diff^{b}(J)=\calo_V\) and
\(\Diff^{b-1}(J)\) is a proper sheaf of ideals).

The closed set of interest is the set of \(b\)-fold points of \(X\)
(i.e. \(V(\Diff^{b-1}(J))\)).  Consider now a \(b\)-permissible
transformation, say \( V'\to V \) as in \ref{eq1}.  Recall that \(
J\calo_{W_1}=I(H')^b J' \).

In this case \(J'\) is the sheaf of ideals defining a hypersurface
\(X'\subset V'\), which is the strict transform of the hypersurface
\(X\).

It is not hard to check that \(J'\) has at most order \(b\) at points
of \(V'\) (i.e. that \(V(\Diff^{b}(J'))=\emptyset)\).  If, in
addition, \(J'\) has no point of order \(b\), then we say that \(\pi\)
defines a {\em \(b\)-simplification} of \(J\).  At any rate, the
closed set of interest is the set of \(b\)-fold points \(X'\).

If \(V(\Diff^{b-1}(J'))\neq \emptyset\), let
\(V'\stackrel{\pi'}{\longleftarrow} V''\) denote the monoidal
transformation with some center \(Y'\subset V(\Diff^{b-1}(J'))\).  So
\(\pi'\) is \(b\)-permissible, and set
\begin{equation*}
J'\calo_{V'}=I(H'')^b J''
\end{equation*}
where \( I(H'') \) is the sheaf of functions vanishing along the
exceptional hypersurface \( H'' \) of \( \pi'' \).

So again \(J''\) has at most points of order \(b\), and if it does,
define a \(b\)-permissible transformation at some smooth center
\(Y''\subset V(\Diff^{b-1}(J'')))\).

So for \(J\) and \(b\) as before, we define, by iteration, a
\(b\)-permissible sequence
\begin{equation} \label{EqSeqPerm}
V\stackrel{\pi }{\longleftarrow} V' \stackrel{\pi'
}{\longleftarrow}V'' \stackrel{\pi'' }{\longleftarrow}\cdots
\stackrel{\pi^{(r-1)}}{\longleftarrow}
V^{(r)}\stackrel{\pi^{(r)}}{\longleftarrow} V^{(r+1)}
\end{equation}
and factorizations \(J^{(i-1)}\calo_{V^{(i)}}=I(H^{(i)})^b J^{(i)}\).
Where \(J^{(i)}\) is the sheaf of ideals defining a hypersurface
\(X^{(i)}\subset V^{(i)}\), which is the strict transform of \(X\).

From the point of view of resolution it is clear that our interest is
to define a \(b\)-permissible sequence so that \(X_{r+1}\) has no
\(b\)-fold points.

We say that a \(b\)-permissible sequence (\ref{EqSeqPerm}) defines a
{\em \(b\)-simplification} of \(J\subset \calo_W\) if the jacobian of
\(V\leftarrow V^{(r+1)}\) has normal crossings, and
\(V(\Diff^{b-1}(J^{(r+1)}))=\emptyset\) (i.e. if \(X^{(r+1)}\) has at
most points of multiplicity \(b-1\)).

Hironaka attaches to the original data \(J\) and \(b\) the pair
\((J,b)\) (\ref{parA}).  The closed set assigned to this pair in \(V\)
is \(\Sing(J,b)=V(\Diff^{b-1}(J))\).  In our case, the \(b\)-fold
points of the hypersurface \(X\).  \medskip

Here we attach to the original data a Rees algebra (up to integral
closure), namely \(\mathcal{G}=\calo_V[JW^b]\).  And to this Rees
algebra a closed set in \(V\), namely \(\Sing(\mathcal{G})\), which is
again \(V(\Diff^{b-1}(J))\).

Moreover, we extended \(\mathcal{G}\) to a Diff-algebra
\(G(\mathcal{G})\), and
\(\Sing(\mathcal{G})=\Sing(G(\mathcal{G}))\) (\ref{thopG}).

Let us focus on one \(b\)-permissible transformation \(\pi\)
(\ref{eq1}).  The transform of Hironaka's pair is the pair \((J',b)\).
The transformation \(\pi\) is also permissible for both
\(\mathcal{G}\) and \(G(\mathcal{G})\), defining transforms of Rees
algebras, say \(\mathcal{G}'\) and \(G(\mathcal{G})'\) on \(V'\).

Note that, in our setting, \(J'\) is the ideal defining  \(X'\),
which is the strict transform of \(X\).  The closed set assigned
to \((J',b)\) is the set of \(b\)-fold points of \(X'\). On the
other hand, \(\mathcal{G}'=\calo_{V'}[J'W^b]\), is such that
\(\Sing(\mathcal{G}')\) is again the set of \(b\)-fold points
\(X'\). A similar relation holds between pairs \((J^{(i)},b)\) and
the Rees algebras \(\mathcal{G}^{(i)}\) (transform of
\(\mathcal{G}\)), for every  \(b\)-permissible sequence
(\ref{EqSeqPerm}).

The natural question is on how do the successive transforms of
\(G(\mathcal{G})\) relate to the transforms of \(\mathcal{G}\).  The
following theorem will address this question.  It proves that the
\(G\)-operator on Rees algebras is, in a natural way, compatible with
transformation.
\begin{theorem} \label{th41} (Giraud)
Let \(\mathcal{G}\) be a Rees algebra on a smooth scheme \(V\), and
let \(V\longleftarrow V'\) be a permissible (monoidal) transformation
for \(\mathcal{G}\).  Let \(\mathcal{G}'\) and \(G(\mathcal{G})'\)
denote the transforms of \(\mathcal{G}\) and \(G(\mathcal{G})\).
Then:
\begin{equation*}
    \mathcal{G}'\subset G(\mathcal{G})'\subset G(\mathcal{G}')
\end{equation*}
so that we have \(G(\mathcal{G}')=G(G(\mathcal{G})')\).
\end{theorem}

\begin{proof}
It is enough to prove that \( G(\mathcal{G})'\subset G(\mathcal{G}') \).

Assume that \( \mathcal{G}=\bigoplus_{n}J_{n}W^{n} \) and
that the monoidal transformation  has center \(
Y\subset\Sing(\mathcal{G}) \).

Let \( x'\in V' \) a point mapped to \( x\in V \).  If \( f\in J_{n}
\) at \( x \) then by hypothesis \( f\in I(Y)^{n} \).  The total
transform of \( f\in\calo_{V',x'} \) is \( f=(h')^{n}f' \) where \( h'
\) is a generator of \( I(H')_{x'} \).  It can be checked that
\begin{equation}
    \Diff_{\calo_{V,x}/k}^{r}(f)\subset
    I(H')^{n-r}_{x'}\Diff_{\calo_{V',x'}/k}^{r}(f')
\end{equation}
Then for every \( D\in\Diff_{\calo_{V,x}/k}^{r} \), we have that for the
total transform \( D(f) \) in \( \calo_{V',x'} \):
\( D(f)=(h')^{n-r}D'(f') \) where \( D'\in\Diff_{\calo_{V',x'}/k}^{r} \).

Assume that \( \mathcal{F}=\{f_{1}W^{n_{1}},\ldots,f_{s}W^{n_{s}}\}
\) is a set of generators of \( \mathcal{G}_{x} \).
Note that the set \( \mathcal{F}'=\{f'_{1}W^{n_{1}},\ldots,f'_{s}W^{n_{s}}\}
\) generates \( \mathcal{G}'_{x'} \).
Recall that a set of generators of \( G(\mathcal{G}) \) is
\begin{equation*}
    \{D(f_{i})W^{m-r}\mid D\in\Diff^{r}_{\calo_{V,x}},\
    0\leq r<m\leq n_{i},\ i=1,\ldots,s\}
\end{equation*}
and a set of generators of the transform \( G(\mathcal{G})' \) is
\begin{equation*}
    \{(h')^{n_{i}-r}D'(f'_{i})W^{m-r}\mid
    0\leq r<m\leq n_{i},\ i=1,\ldots,s\}
\end{equation*}
where \( D' \) are some differential operators in \(
\Diff_{\calo_{V',x'}}^{r} \).
\end{proof}

\begin{definition} \label{DefResolRA}
A resolution of a Rees algebra \( \mathcal{G} \) is a sequence of
transformations, say
\begin{equation}
    \begin{array}{ccccccccc}
        V & \leftarrow & V' & \leftarrow & \cdots & \leftarrow &
    V^{(r)} & \leftarrow & V^{(r+1)}  \\
        \mathcal{G} &  & \mathcal{G}' &  &  &  & \mathcal{G}^{(r)} &
    & \mathcal{G}^{(r+1)}
    \end{array}
    \label{EqSeqRees}
\end{equation}
such that \( \Sing(\mathcal{G}^{(r+1)})=\emptyset \), and the
exceptional locus of $V\leftarrow  V^{(r+1)}$ is a union of smooth
hypersurfaces with normal crossings.
\end{definition}

\section{Idealistic exponents versus basic objects.}
\label{sec5}

Recall the definition of idealistic equivalence \ref{defe}.
\begin{proposition} \label{PropIdEquiv}
Let \((J_{1},b_{1})\) and \((J_{2},b_{2})\) be
idealistic equivalent.  Then:
\begin{enumerate}
    \item \(\Sing(J_{1},b_{1})=\Sing(J_{2},b_{2})\).

    \noindent Note, in particular, that every monoidal transform
    \(V\leftarrow V'\) on a center \(Y\subset
    \Sing(J_{1},b_{1})=\Sing(J_{2},b_{2})\) defines transforms, say
    \((J'_1,b_{1})\) and \((J'_2,b_{2})\) on \(V'\).

    \item The pairs \((J'_1,b_{1})\) and \((J'_2,b_{2})\) are
    idealistic equivalent on \(V'\).
\end{enumerate}
\end{proposition}

\begin{parrafo}
If two pairs \((J_{1},b_{1})\) and \((J_{2},b_{2})\) are
idealistic equivalent over \(V\), the same holds for the
restrictions to every open subset of \(V\), and also for
restrictions in the sense of \'etale topology, and even for smooth
topology (i.e. pull-backs by smooth morphisms \(V'\to V\)).

Note that if \((J_{1},b_{1})\) and \((J_{2},b_{2})\) are
idealistic equivalent, they define the same closed set on \(V\)
(i.e. \(\Sing(J_{1},b_{1})=\Sing(J_{2},b_{2})\)), and the same
holds for monoidal transformations, pull-backs by smooth schemes,
and hence by concatenation of both kinds of transformations.  When
this last condition holds on the singular locus of two pairs we
say that they {\em define the same closed sets}.
\end{parrafo}

\begin{definition} \label{DefEquivBas}
Two pairs \((J_{1},b_{1})\) and \((J_{2},b_{2})\) are {\em
basically} equivalent on \(V\), if they define the same closed
sets.
\end{definition}
Proposition~\ref{PropIdEquiv} says that if two pairs are idealistic
equivalent over \(V\), then they are basically equivalent.  \medskip

An {\em idealistic exponent}, as defined by Hironaka in
\cite{Hironaka1977}, is an equivalence class of pairs in the sense of
idealistic equivalence.  Whereas the notion of equivalence among basic
objects (see \cite{Villamayor1989} or \cite{Villamayor1992}) is
\ref{DefEquivBas}.  In fact, the key point for constructive
desingularization was to define an algorithm of resolutions of pairs
\((J,b)\), so that two basically equivalent pairs undergo exactly the
same resolution.  \medskip

Recall now the definition of integrally equivalence on Rees algebras
\ref{defei}.
\begin{proposition}
Let \(\mathcal{G}_{1}\) and \(\mathcal{G}_{2}\) be two integrally
equivalent Rees algebras over \(V\). Then:
\begin{enumerate}
    \item \(\Sing(\mathcal{G}_{1})=\Sing(\mathcal{G}_{2})\).

    \noindent Note, in particular, that every monoidal transform
    \(V\leftarrow V'\) on a center
    \(Y\subset\Sing(\mathcal{G}_{1})=\Sing(\mathcal{G}_{2})\) defines
    transforms, say \(\mathcal{G}'_1\) and \(\mathcal{G}'_2\) on
    \(V'\).

    \item \(\mathcal{G}'_1\) and \(\mathcal{G}'_2\) are integrally
    equivalent on \(V'\).
\end{enumerate}
\end{proposition}

If \(\mathcal{G}_{1}\) and \(\mathcal{G}_{2}\) are {\em integrally}
equivalent on \(V\), the same holds for every open restriction, and also
for pull-backs by smooth morphisms \(V'\to V\).

On the other hand, as \(\mathcal{G}'_1\) and \(\mathcal{G}'_2\) are
integrally equivalent,  they define the same closed set on \(V'\)
(the same singular locus), and the same holds for further monoidal
transformations, pull-backs by smooth schemes, and concatenations of
both kinds of transformations.

When this condition holds on the singular locus of two Rees
algebras over \(V\), we say that they {\em define the same closed
sets}.

\begin{definition}
Two Rees algebras over \(V\), say \(\mathcal{G}_{1}=\bigoplus_{n\geq
0}I_nW^n\) and \(\mathcal{G}_{2}=\bigoplus_{n\geq 0}J_nW^n\), are
\emph{basically} equivalent, if both define the same closed sets.
\end{definition}

The previous Proposition asserts that if
\(\mathcal{G}_{1}=\bigoplus_{n\geq 0}I_nW^n\) and
\(\mathcal{G}_{2}=\bigoplus_{n\geq 0}J_nW^n\) are integrally
equivalent, then they are basically equivalent.

\begin{proposition}
Let \( (J_{1},b_{1}) \) and \( (J_{2},b_{2}) \) be two pairs.
\begin{enumerate}
    \item The pairs \((J_{1},b_{1})\) and \((J_{2},b_{2})\) are
    idealistically equivalent over a smooth scheme \(V\), if and only
    if the Rees algebras \(\mathcal{G}_{(J_{1},b_{1})}\) and
    \(\mathcal{G}_{(J_{2},b_{2})}\) are integrally equivalent
    (\ref{parC} ).

    \item The pairs \((J_{1},b_{1})\) and \((J_{2},b_{2})\) are
    basically equivalent over \(V\), if and only if the Rees algebras
    \(\mathcal{G}_{(J_{1},b_{1})}\) and
    \(\mathcal{G}_{(J_{2},b_{2})}\) are basically equivalent.
\end{enumerate}
\end{proposition}

\begin{proof}
Note that (1) was already proved in \ref{PropEquivGJ}.  (2) follows
from the fact that if \( (J',b) \) is the transform of \( (J,b) \)
then \( \mathcal{G}_{(J',b)} \) is the transform of \(
\mathcal{G}_{(J,b)} \).
\end{proof}

\section{Functions on pairs and simple Rees algebras}

\begin{definition} \label{Defusc}
Let \( V \) be a noetherian topological space
\cite[p.~5]{Hartshorne1977} and let \( (\Lambda,\leq) \) be a
totally ordered set.  A function \(
f:V\longrightarrow(\Lambda,\leq) \) is said to be
upper-semi-continuous (u.~s.~c.)  if it takes finitely many
values, and for every \( \alpha\in\Lambda \) the set \( \{x\in V\mid
f(x)\geq\alpha\} \) is closed.  \medskip

We denote by \( \max{f} \) the biggest value and \( \MaxB{f}=\{x\in
V\mid f(x)=\max{f}\} \), which is closed in \( V \).
\end{definition}

\begin{remark}\label{meolv}
The goal is to define an u.s.c function, say \( f_{\mathcal{G}}
\), for every Rees algebra \( \mathcal{G} \),  such that \(
\MaxB{f_{\mathcal{G}}} \) is a permissible center. Set \( V'\to V
\) the transformation with center \( \MaxB{f_{\mathcal{G}}} \).
Then the function \( f_{\mathcal{G}'} \) is such that \(
\max{f_{\mathcal{G}}}>\max{f_{\mathcal{G}'}} \).

An algorithm of resolution of Rees algebras is the assignment of
functions \( f_{\mathcal{G}} \), for each \( \mathcal{G} \), such
that a sequence (\ref{EqSeqRees}) is defined by setting each
transformation in \ref{EqSeqRees} with center \(
\MaxB{f_{\mathcal{G}^{(i)}}} \). Furthermore, inequalities \(
\max{f_{\mathcal{G}}}>\max{f_{\mathcal{G}'}}>\cdots
>\max{f_{\mathcal{G}^{(r)}}} \) will hold,
and the sequence (\ref{EqSeqRees}) defined
in this way is a
resolution for suitable \( r=r_{\mathcal{G}} \).
\end{remark}

\begin{parrafo}\label{agre1}
We introduce a function, again a natural analog to that defined
for idealistic exponents.  Fix \(x\in V\).  Given \(fW^n\in
I_nW^n\), where \(\mathcal{G}=\bigoplus_{n\geq 0}I_nW^n\), set
\begin{equation*}
\ord_x(fW^{n})=\frac{\nu_x(f)}{n}\in \mathbb{Q};
\end{equation*}
called the order of \(f\) (weighted by \(n\)), where \(\nu_x\)
denotes the order at the local regular ring \(\calo_{V,x}\).  Note
that \(x\in \Sing(\mathcal{G})\) if and only if \(\ord_x
(fW^{n})\geq 1\) for all \( n\geq 1 \). We also define
\begin{equation*}
\ord_x (\mathcal{G})=
\inf\{\ord_x (fW^{n})\mid fW^n\in I_nW^n, \quad n\geq 1\}.
\end{equation*}

So, in general \(\ord_x (\mathcal{G})\geq 1\) iff \(x\in
\Sing(\mathcal{G})\).  A Rees algebra \(\mathcal{G}\) is said to
be simple at \(x\) if \(\ord_x (\mathcal{G})= 1\).
\end{parrafo}

\begin{proposition}\label{3propsingZ}
\begin{enumerate}
    \item If \(\mathcal{G}\) is a Rees algebra generated over
    \(\calo_V\) by
    \(\mathcal{F}=\{ g_{i}W^{n_i}\mid 1\leq i\leq m \}\), then
    \( \ord_x (\mathcal{G})=
    \min\{\ord_x (g_{i}W^{n_{i}})\mid 1\leq i\leq m\} \).
    If \(N\) is a common multiple of all \(n_i\), \(1\leq i\leq
    m\), then
    \begin{equation*}
    \ord_x (\mathcal{G})=\frac{\nu_x(I_N)}{N}.
    \end{equation*}

    \item If \(\mathcal{G}_{1}\) and \( \mathcal{G}_{2}\) are Rees algebras with the same integral closure (e.g. if
    \(\mathcal{G}_{1}\subset \mathcal{G}_{2}\) is a finite extension),
    then, for all  \(x\in
    \Sing(\mathcal{G}_{1})(=\Sing(\mathcal{G}_{2}))\)
    \begin{equation*}
    \ord_x (\mathcal{G}_{1})=\ord_x (\mathcal{G}_{2}).
    \end{equation*}
    In particular, the function is compatible with integral
    equivalence (see \ref{Require}).

    \item Let \(G(\mathcal{G})\) be the extension of \(\mathcal{G}\)
    to a differential Rees algebra relative to \(k\), then for all  \(x\in
    \Sing(\mathcal{G})(=\Sing(G(\mathcal{G})))\).
    \begin{equation*}
    \ord_x (\mathcal{G})=\ord_x (G(\mathcal{G})).
    \end{equation*}
\end{enumerate}
\end{proposition}

\begin{remark}
The function \(\ord\) was introduced by Hironaka in the context of
pairs.  Given a pair \((J,b)\) as in \ref{parA}; and assume that
\(J\subset \calo_V\) is a non-zero sheaf of ideals, a function \(\ord:
\Sing(J,b) \to \mathbb{Q}\) is defined by setting
\begin{equation*}
\ord_x(J,b)= \frac{\nu_{x}(J)}{b}.
\end{equation*}

Note that if \(\mathcal{G}_{(J,b)}\) is the Rees algebra attached to
\((J,b)\), then \(\Sing(J,b)=\Sing(\mathcal{G}_{(J,b)})\), and for all
\(x\in \Sing(J,b)\): \( \ord_x(J,b)=\ord_x(\mathcal{G}_{(J,b)})\)
(\ref{parC}).
\end{remark}

\begin{proposition} \label{PropOrdBaja}
Let \( \mathcal{G} \) be a Rees algebra and \( V'\to V \) be a
permissible transformation with center \( Y \). Denote by \(
\mathcal{G}' \) the transform of \( \mathcal{G} \).
\begin{enumerate}
    \item If the function \( \ord(\mathcal{G}) \) is constant along
    \( Y \) then, for all \( x'\in V' \) mapping to \( x\in V \):
    \begin{equation*}
        \ord_{x}(\mathcal{G})\geq\ord_{x'}(\mathcal{G}')
    \end{equation*}

    \item  If \( \mathcal{G} \) is simple at $x$, then \( \mathcal{G}' \) is
    simple at $x'$.

    \item  If \( \codim_{V}(Y)=1 \) then \( V'\to V \) is an
    isomorphism but \( \mathcal{G}' \) and \( \mathcal{G} \) are
    different, in fact
    \begin{equation*}
        \ord_{x'}(\mathcal{G})=\ord_{x}(\mathcal{G})-1, \qquad
    \forall x=x'\in V=V'
    \end{equation*}
\end{enumerate}
\end{proposition}

\begin{proof}
If \( \mathcal{G}=\mathcal{G}_{(J,b)} \) then \(
\mathcal{G}'=\mathcal{G}_{(J',b)} \) and it is well known that \(
\nu_{x}(J)\geq\nu_{x'}(J') \).

In the general case \( \mathcal{G} \) is integral over some \(
\mathcal{G}_{(J_{N},N)} \) (\ref{rkK1}), and (1) follows from
\ref{3propsingZ}.

(2) is a consequence of (1), and (3) follows from the definition
of transformation \ref{parwt}.
\end{proof}

\begin{proposition}
Assume that the ground field \( k \) is of characteristic zero.  Fix
\( x\in \Sing(\mathcal{G}) (\subset V) \).  If \(
\ord(\mathcal{G})(x)=1 \) then there is, locally at \( x \), a smooth
hypersurface \( Z \) such that
\begin{enumerate}
    \item \( \calo_{V}[I(Z)W]\subset G(\mathcal{G}) \).  And hence \(
    \Sing(\mathcal{G})\subset Z \).

    \item If \( V'\to V \) is a blow-up at a permissible center \(
    C\subset\Sing(\mathcal{G}) \), and if \(
    \mathcal{G}' \) is the transform of \( \mathcal{G} \) and \(
    Z'\subset V' \) is the strict transform of \( Z \), then
    \begin{equation*}
    \calo_{V'}[I(Z')W]\subset G(\mathcal{G}').
    \end{equation*}
\end{enumerate}
\end{proposition}
\proof Fix a regular system of parameters \(\{x_1,\dots,x_d\}\) at
\(\calo_{V,x}\), and an element \(f_nW^n\in \mathcal{G}\) so that
\(f_n\) has order \(n\) at \(\calo_{V,x}\).  Note that
\(\Delta^{\alpha}(f_n)W\in G(\mathcal{G})\) for \(|\alpha|=n-1\),
which can be chosen so that \(\Delta^{\alpha}(f_n)\) has order one,
defining a smooth hypersurface \(Z\) locally at \(x\) (\ref{parod}).

\begin{theorem}
\textbf{\label{ThJarek} }(W{\l}odarczyk, \cite{Wlodarczyk2005})
Let \( x\in V \) be a simple point of \(\mathcal{G}\), and assume
that locally at \( x \), there are two hypersurfaces \(
Z_{1},Z_{2}\subset V \) such that
\begin{equation*}
    \calo_{V}[I(Z_{i})W]\subset G(\mathcal{G}) \qquad i=1,2
\end{equation*}
Then there exist  \'etale neighborhoods $ \phi_1$, $\phi_2$ : $U\to V$ of \( x =\phi_1(y)=\phi_2(y)\), where $y \in U$,  such that \( \phi^*_1(Z_{1})= \phi^*_2(Z_{2}) \)
and \( \phi^{\ast}_1(G(\mathcal{G}_{}))=\phi^{\ast}_2(G(\mathcal{G}_{}))\).
\end{theorem}

\begin{parrafo}
If \( \mathcal{G} \) is simple at \( x\in\Sing(\mathcal{G}) \),
then resolution of \( \mathcal{G} \) will reduce, locally, to the
resolution of the restriction, say \( G(\mathcal{G})_{Z} \), of \(
G(\mathcal{G}) \) to a hypersurface \( Z \). Where $Z$ is such
that \( \calo_{V}[I(Z)W]\subset G(\mathcal{G}) \).
By
theorem~\ref{ThGEntera} this procedure is well defined up to
integral closure.  If \( \mathcal{G}_{1} \) and \( \mathcal{G}_{2}
\) have the same integral closure then \( G(\mathcal{G}_{1})_{Z}
\) and \( G(\mathcal{G}_{2})_{Z} \) have the same integral
closure.

By theorem~\ref{ThJarek} this procedure does not depend on the choice
of the hypersurface \( Z \).  In fact, if two possible hypersurfaces
\( Z_{1} \) and \( Z_{2} \) fulfill the previous condition, then the
restrictions \( G(\mathcal{G})_{Z_{1}} \) and \(
G(\mathcal{G})_{Z_{2}} \) are the same (\'etale locally at \( y \)).
\end{parrafo}

\section{Reduction to the simple case}

Given a Rees algebra $\mathcal{G}$ we have defined \(
\ord_{x}(\mathcal{G})\in \mathbb{Q}\geq 1 \) for every point \(
x\in\Sing(\mathcal{G}) \) . Suppose that \(
\mathcal{G}=\mathcal{G}_{(J,b)} \) and set \(
\omega=\ord_{x}(J,b)\geq 1 \).  Given an idealistic pair \( (J,b)
\) we may consider a new pair locally at \( x \), with order one,
say \( (J,a) \), where \( a=\omega b \).  The next definition is
the analogous formulation, now for the case of Rees algebras.
\begin{definition} \label{DefTwist}
Let \( \mathcal{G}=\bigoplus_{n\geq 0}J_{n}W^{n} \) be a Rees algebra
on \( V \) and fix \( \omega\in\mathbb{Q} \), \( \omega>0 \).  We
define the twisted algebra \( \mathcal{G}(\omega) \) as follows
\begin{equation*}
    \mathcal{G}(\omega)= \bigoplus_{n\geq 0}J_{\frac{n}{\omega}}W^{n}
\end{equation*}
where we set \( J_{\frac{n}{\omega}}=0 \) whenever \(
\dfrac{n}{\omega} \) is not an integer.

More precisely: if \( \omega=\dfrac{a}{b} \) where \( a \) and \( b \)
are integers with no common factors, and if \( n=am \) for some \(
m\in\mathbb{Z} \), then \( J_{\frac{n}{\omega}}=J_{bm} \) and \(
J_{\frac{n}{\omega}}=0 \) if \( a \) not divide \( n \).
\end{definition}

\begin{remark} \label{RemTwist}
 \( \mathcal{G}(\omega) \) can be constructed as follows: If \( \omega=\dfrac{a}{b}
\) with \( (a,b)=1 \), then consider the isomorphism \(
\Phi:\calo_{V}[W^{b}]\to\calo_{V}[W^{a}] \) by sending \( W^{b} \)
to \( W^{a} \).  Then the twisted Rees algebra is
\begin{equation*}
    \mathcal{G}(\omega)=
    \Phi\left(V^{(b)}\left(\mathcal{G}\right)\right)
\end{equation*}
\end{remark}

\begin{proposition}\label{rk73}
The definition~\ref{DefTwist} fulfills the requirement of
(\ref{Require}).
\begin{enumerate}
    \item If \(\mathcal{G}=\mathcal{G}_{(J,b)}\) and \(
    \omega\in\mathbb{Q} \) is such that \( b\omega\in\mathbb{Z} \)
    then \( \mathcal{G}(\omega)=\mathcal{G}_{(J,b\omega)}\).

    \item If \( \mathcal{G}_{1} \) and \( \mathcal{G}_{2} \) have the
    same integral closure then \( \mathcal{G}_{1}(\omega) \) and \(
    \mathcal{G}_{2}(\omega) \) have the same integral closure.
\end{enumerate}
\end{proposition}

\begin{proof}
(1) follows straightforward.  (2) follows from \ref{RemTwist} and
\ref{paramigo}, which asserts that the operator $V^m$ is
compatible with integral closure.
\end{proof}

\begin{proposition}
\( \mathcal{G}(\omega) \) is a Rees algebra, and \(
\omega\ord_{x}(\mathcal{G}(\omega))=\ord_{x}(\mathcal{G}) \) at
every point \( x\in V \).
\end{proposition}

\begin{proof}
Set $\mathcal{G}=\bigoplus I_nW^n$. As $\mathcal{G}$ is finitely
generated it is a finite extension of $\mathcal{G}_{(I_N,N)}$ for
a suitable choice of $N$. \ref{rk73} shows that
$\mathcal{G}(\omega)$ is a finite extension of
$\mathcal{G}_{(I_N,\omega N)}$, therefore \( \mathcal{G}(\omega)
\) is a Rees algebra.

\medskip

For all  \( x\in V \) we have
\begin{align*}
    \alpha=\ord_{x}(\mathcal{G}) \Longleftrightarrow &
    \nu_{x}(J_{n})\geq \alpha n \ \forall n>1, \text{ and } \exists\
    n_{0} \text{ with } \nu_{x}(J_{n_{0}})=\alpha n_{0} \\
    \Longleftrightarrow & \nu_{x}(J_{bn})\geq \alpha bn \ \forall n>1
    \text{ and } \exists\ n_{0} \text{ with }
    \nu_{x}(J_{bn_{0}})=\alpha bn_{0} \\
    \Longleftrightarrow & \nu_{x}(J_{\frac{an}{\omega}})\geq
    \frac{\alpha an}{\omega} \text{ and } \exists\ n_{0} \text{ with }
    \nu_{x}(J_{\frac{an_{0}}{\omega}})=\frac{\alpha an_{0}}{\omega} \\
    \Longleftrightarrow &
    \frac{\alpha}{\omega}=\ord(\mathcal{G}(\omega))
\end{align*}
\end{proof}

Note also that if \( \omega\geq 1 \) then \(
\Sing(\mathcal{G}(\omega))\subset\Sing(\mathcal{G}) \), but if \(
\omega<1 \) then this inclusion may not be satisfied.

\begin{corollary}
If \( \omega=\ord_{x}(\mathcal{G}) \) then \(
\ord_{x}(\mathcal{G}(\omega))=1 \).
\end{corollary}

In order to achieve a resolution of a Rees algebra \( \mathcal{G} \),
we will attach to \( \mathcal{G} \) some additional data.  The first
one will be a set of irreducible hypersurfaces having normal
crossings, say \( E=\{H_{1},\ldots,H_{r}\} \).
\begin{definition} \label{DefDivMonom}
Let \( \mathcal{G}=\bigoplus_{n}J_{n}W^{n} \) be a Rees algebra,
\( E=\{H_{1},\ldots,H_{r}\} \) a set of irreducible hypersurfaces
having normal crossings, and \(
\underline{a}=(a_{1},\ldots,a_{r})\in\mathbb{Q}^{r} \).

We say that the formal monomial \( I(H_{1})^{a_{1}}\cdots
I(H_{r})^{a_{r}} \) divides the algebra \( \mathcal{G} \) (or \(
E^{\underline{a}} \) divides \( \mathcal{G} \)) if, for every \( n
\):
\begin{equation*}
    J_{n}\subset I(H_{1})^{\lfloor na_{1}\rfloor}\cdots
    I(H_{r})^{\lfloor na_{r}\rfloor}
\end{equation*}
where \( \lfloor na_{i}\rfloor\leq na_{i} \) is the integer part of \(
na_{i} \), \( i=1,\ldots,r \).
\end{definition}

\begin{parrafo} \label{DivRequire}
We now show that Definition~\ref{DefDivMonom} fulfills the requirement
of \ref{Require}.
\begin{enumerate}
    \item Let \( \mathcal{G}=\mathcal{G}_{(J,b)} \) and \(
    \underline{a}\in\mathbb{Q}^{r} \) be such that \(
    b\underline{a}\in\mathbb{Z}^{r} \).
The monomial \( E^{\underline{a}} \) divides \( \mathcal{G} \) if
    and only if the monomial \( E^{b\underline{a}} \) divides the
    ideal \( J \).

    \item Let \( \mathcal{G}_{1} \) and \( \mathcal{G}_{2} \) be two
    Rees algebras which are integrally equivalent.
Then \( E^{\underline{a}} \) divides \( \mathcal{G}_{1} \) if and
    only if \( E^{\underline{a}} \) divides \( \mathcal{G}_{2} \).
\end{enumerate}
\end{parrafo}

\begin{parrafo} \label{MaximalQ}
Let \( N \) be such that \(\mathcal{G}_{(J_{N},N)}\subset
\mathcal{G}=\bigoplus_{n}J_{n}W^{n}\) is finite. Set
\(\alpha_1,\dots,\alpha_r\) be the biggest integers for which
\(J_{N}\subset I(H_{1})^{\alpha_1}\cdots I(H_{r})^{\alpha_r}\) and set
\(a_i=\dfrac{\alpha_i}{N}\in \mathbb{Q}\).  Then \(
\underline{a}=(a_{1},\ldots,a_{r})\in\mathbb{Q}^{r} \) is the maximal
vector such that the monomial \( E^{\underline{a}} \) divides \(
\mathcal{G} \). In the sense that if \(
\underline{a}'\neq\underline{a} \) and \( a'_{i}\geq a_{i} \) then \(
E^{\underline{a}'} \) does not divide \( \mathcal{G} \).

We claim that this assertion is independent of the choice of \(N\)
with the previous property.  Note that \(\alpha_i\) is the valuation
of \(J_N\) at the generic point of the hypersurface.  On the other
hand Proposition  \ref{PropEquivGJ} says that if a Rees algebra
\(\mathcal{G}\) is integrally equivalent to two different algebras
attached to two pairs, say \(\mathcal{G}_{(J_{1},b_{1})}\) and to
\(\mathcal{G}_{(J_{2},b_{2})}\), then the pairs \((J_{1},b_{1})\) and
\((J_{2},b_{2})\) are idealistic equivalent (\ref{defe}).  In
particular \(J_1^{b_2}\) and \(J_2^{b_1}\) vanish along \(H_i\) with
the same order.
\end{parrafo}

\begin{definition} \label{DefWeak}
Let \( \mathcal{G}=\bigoplus_{n}J_{n}W^{n} \) be a Rees algebra
and \( E \) be a set of irreducible hypersurfaces having only
normal crossings. Let \( \underline{a} \) be the maximal vector
such that \( E^{\underline{a}} \) divides \( \mathcal{G} \). For
all \( n \) there is an ideal \( I_{n} \) such that
\begin{equation*}
    J_{n}=I(H_{1})^{\lfloor na_{1}\rfloor}\cdots I(H_{r})^{\lfloor
    na_{r}\rfloor} I_{n}
\end{equation*}
Set \( \mathcal{G}^{\weak}=\bigoplus_{n}I_{n}W^{n} \) and define the
function \( \word(\mathcal{G}):V\to\mathbb{Q} \) to be \(
\word(\mathcal{G})=\ord(\mathcal{G}^{\weak}) \).
\end{definition}

It follows from \ref{MaximalQ} that \( \mathcal{G}^{\weak} \) is
finitely generated, and hence a Rees algebra.  Using notation as
in \ref{MaximalQ}, we have that \( J_{N}=I(H_{1})^{\alpha_1}\cdots
I(H_{r})^{\alpha_r}I_{N} \) and \(
\mathcal{G}_{(I_{N},N)}\subset\mathcal{G}^{\weak} \) is integral.

\begin{parrafo} 
Definition \ref{DefWeak} satisfies requirement in \ref{Require}.
\begin{enumerate}
    \item Let \( (J,b) \) be a pair such that
    \begin{equation*}
    J=I(H_{1})^{\alpha_{1}}\cdots I(H_{r})^{\alpha_{r}}I
    \end{equation*}
    and \( I\not\subset I(H_{i}) \) for \( i=1,\ldots,r \).  If \(
    \mathcal{G}=\mathcal{G}_{(J,b)} \) then \(
    \mathcal{G}^{\weak}=\mathcal{G}_{(I,b)} \).

    Moreover, this function \( \word(\mathcal{G}) \) coincides with the function \(
    \word_{(J,b)} \), defined in terms of the pair \((J,b)\) in
    \cite{BravoEncinasVillamayor2005}.

    \item If \( \mathcal{G}_{1} \) and \( \mathcal{G}_{2} \) are
    integrally equivalent then \( \mathcal{G}_{1}^{\weak} \) and \(
    \mathcal{G}_{2}^{\weak} \) are integrally equivalent.

    In particular \( \word(\mathcal{G}_{1})=\word(\mathcal{G}_{2})
    \) (both functions are equal).
\end{enumerate}
\end{parrafo}

\begin{parrafo} \label{TransWeak}
Let \( V'\to V \) be a transformation with irreducible center \( C
\). We assume that \( C\subset\Sing(\mathcal{G}) \), and that the
function \( \word(\mathcal{G}) \) is constant along \( C \) and
takes the value, say \( \omega\in\mathbb{Q} \).

Let \( \underline{a} \) be the maximal vector such that \(
E^{\underline{a}} \) divides \( \mathcal{G} \).  Let \(
\mathcal{G}' \) be the transform of \( \mathcal{G} \), and set \(
E'=\{H'_{1},\ldots,H'_{r},H'_{r+1}\} \) where \( H'_{r+1} \) is
the exceptional divisor of \( V'\to V \), and \( H'_{i} \) is the
strict transform of \( H_{i} \), \( i=1,\ldots,r \).  Then:

\begin{enumerate}
    \item The monomial \( (E')^{\underline{a}'} \) is the maximal
    monomial which divides \( \mathcal{G}' \).

    Where \( \underline{a}'=(a'_{1},\ldots,a'_{r},a'_{r+1}) \), \(
    a'_{i}=a_{i} \) for \( i=1,\ldots,r \) and
    \begin{equation*}
    a'_{r+1}=\omega-1+\sum_{H_{i}\supset C}a_{i}
    \end{equation*}

    \item The twisted Rees algebra \( \mathcal{G}^{\weak}(\omega) \)
    is simple and its transform is
    \begin{equation*}
    \left(\mathcal{G}^{\weak}(\omega)\right)'=
    \left(\mathcal{G}'\right)^{\weak}(\omega)
    \end{equation*}
\end{enumerate}
\end{parrafo}

\begin{parrafo}
We add some more information to a the Rees algebra in addition to
\( E \), say \( (\mathcal{G},D,E) \) where \( D\subset E \) is a
suitable subset of hypersurfaces. Given a monoidal transformation
we define the transform of \( (\mathcal{G},D,E) \) to be \(
(\mathcal{G}',D',E') \) where \( \mathcal{G}' \) is the transform
of \( \mathcal{G} \), \( E' \) is as in \ref{TransWeak} and \( D'
\) consist of the strict transforms of \( D \) if \(
\max\word(\mathcal{G})=\max\word(\mathcal{G}') \) and \( D'=E' \)
otherwise.
\end{parrafo}

\begin{definition}
Given \( (\mathcal{G},D,E) \) we define the function \(
t(\mathcal{G})=(\word(\mathcal{G}),n(\mathcal{G})) \) where \(
n(\mathcal{G})(x) \) is the number of hypersurfaces of \( D \) passing
through \( x\in V \).
\end{definition}

\begin{parrafo}
The function \( t(\mathcal{G}) \) satisfies the requirements of
\ref{Require}.
\begin{enumerate}
    \item If \( \mathcal{G}=\mathcal{G}_{(J,b)} \) then \(
    t(\mathcal{G}) \) is the function \( t(J,b) \) defined in terms of
    the pair in \cite{BravoEncinasVillamayor2005}.

    \item If \( \mathcal{G}_{1} \) and \( \mathcal{G}_{2} \) are
    integrally equivalent then \(
    t(\mathcal{G}_{1})=t(\mathcal{G}_{2}) \).

\end{enumerate}
\end{parrafo}

\begin{parrafo}
Given \( (\mathcal{G},D,E) \), if \( \max{t(\mathcal{G})}=(\omega,m)
\) then set:
\begin{equation*}
    \mathcal{T}(\mathcal{G})=
    \mathcal{G}\odot\mathcal{G}^{\weak}(\omega)\odot\mathcal{D}_{m}
\end{equation*}
where \( \mathcal{D}_{m} \) is the Rees ring of the ideal
\begin{equation*}
    \prod_{\{H_{1},\ldots,H_{m}\}\subset D}\left(
    \sum_{j=1}^{m}I(H_{j}) \right)
\end{equation*}
see \cite[definition 15.14]{BravoEncinasVillamayor2005}.

The Rees algebra \( \mathcal{T}(\mathcal{G}) \) is simple and \(
\Sing(\mathcal{T}(\mathcal{G}))=\MaxB{t(\mathcal{G})} \).
\end{parrafo}

\begin{parrafo}
The algebra \( \mathcal{T}(\mathcal{G}) \) satisfies the
requirements of \ref{Require}.

\begin{enumerate}
    \item If \( \mathcal{G}=\mathcal{G}_{(J,b)} \), then \(
    \mathcal{T}(\mathcal{G})=\mathcal{G}_{t(J,b)} \).  Where here \( t(J,b)
    \) denotes the pair attached to \( (J,b) \) in
    \cite[15.14.2]{BravoEncinasVillamayor2005}; defined in terms of the
    function \( t \) of the pair \( (J,b) \).

    \item If \( \mathcal{G}_{1} \) and \( \mathcal{G}_{2} \) are
    integrally equivalent then \( \mathcal{T}(\mathcal{G}_{1}) \) and
    \( \mathcal{T}(\mathcal{G}_{2}) \) are integrally equivalent.
\end{enumerate}
\end{parrafo}

\begin{parrafo}
If \( V'\to V \) is a transformation with center \(
C\subset\Sing(\mathcal{T}(\mathcal{G}))=\MaxB{t(\mathcal{G})} \),
\( (\mathcal{G}',D',E') \) is the transform of \(
(\mathcal{G},D,E) \), and \( \mathcal{T}(\mathcal{G})' \) is the
transform of \( \mathcal{T}(\mathcal{G}) \), then we have the same
relation among the transforms, namely:
\begin{equation*}
    \mathcal{T}(\mathcal{G})'=\mathcal{G}'\odot
    \left(\mathcal{G}'\right)^{\weak}(\omega)\odot\mathcal{D}'_{m}
\end{equation*}
So that
\begin{enumerate}
    \item If \( \max{t(\mathcal{G})}=\max{t(\mathcal{G}')} \) then \(
    \mathcal{T}(\mathcal{G})'=\mathcal{T}(\mathcal{G}') \).

    \item If \( \max{t(\mathcal{G})}>\max{t(\mathcal{G}')} \) then \(
    \Sing(\mathcal{T}(\mathcal{G})')=\emptyset \) and \(
    \mathcal{T}(\mathcal{G}') \) is a new Rees algebra such that \(
    \Sing(\mathcal{T}(\mathcal{G}'))=\MaxB{t(\mathcal{G}')} \).
\end{enumerate}
\end{parrafo}

\begin{theorem}\label{teofin}
The function \( t \) extends naturally, by induction on the
dimension of the ambient space, to an algorithm of resolution of
Rees algebras (\ref{meolv}).

If \( \mathcal{G}_{1} \) and \( \mathcal{G}_{2} \) are integrally
equivalent then the algorithm defines the same resolution for both
Rees algebras.
\end{theorem}

\begin{proof}
The algorithm is defined by induction on \( \dim{V} \).

If \( \dim{V}=1 \), then \( \Sing(\mathcal{G}) \) consists of
finitely many points, and every such point is the center of the
transformation. Theorem follows in this case from
\ref{PropOrdBaja}(3).
\medskip

Assume that an algorithm is defined for Rees algebras over smooth
schemes of dimension \( n-1 \).

Let \( \mathcal{G} \) be a Rees algebra on \( V \), \( \dim{V}=n
\). The Rees algebra \( \mathcal{T}(\mathcal{G}) \) is simple, so
that locally we may choose a smooth hypersurface \( Z\subset V \)
and consider the resolution of the restriction \(
G(\mathcal{T}(\mathcal{G}))_{Z} \) of \(
G(\mathcal{T}(\mathcal{G})) \) to \( Z \), which exists by
induction. W{\l}odarczyks theorem \ref{ThJarek}  asserts that such
local procedures globalize.

In particular there is a sequence of blow-ups, say \( V'\to V \),
such that the transform \( \mathcal{T}(\mathcal{G})' \) has empty
singular locus.  This means that the transform \( \mathcal{G}' \)
is such that \( \max{t(\mathcal{G})}>\max{t(\mathcal{G}')} \).

We continue, similarly, with the simple Rees algebra \(
\mathcal{T}(\mathcal{G}') \), and so on.

This procedure defines a sequence of transformations for \(
\mathcal{G} \), such that the function \( t(\mathcal{G}) \) drops
after finitely many steps.

Now this procedure shall stop since the function \( t(\mathcal{G})
\) may not drop infinitely many times.  This follows from the fact
that \( t(\mathcal{G}) \) is the function associated to some
suitable pair \( (J_{N},N) \); and for pairs the function \( t \)
may not drop infinitely many times, see
\cite{BravoEncinasVillamayor2005}.
\medskip

Last assertion of the theorem follows from the fact that all
constructions satisfy requirements in \ref{Require}.
\end{proof}

\end{document}